\def\softd{{\leavevmode\setbox1=\hbox{d}%
\hbox to 1.05\wd1{d\kern-0.4ex{\char039}\hss}}}
\def\softt{{\leavevmode\setbox1=\hbox{t}%
\hbox to \wd1{t\kern-0.6ex{\char039}\hss}}}
\def\softl{l\kern-0.45ex\raise0.1ex\hbox{'}\kern-0.10ex}
\def\softL{L\kern-0.8ex\raise0.1ex\hbox{'}\kern0.1ex}

\documentstyle[a4,psfig]{article}

\begin{document}
\pagestyle{empty}
\vspace*{40mm}

\begin{center}
{\Large {\bf{MODELLING AND ANALYSIS OF FRACTIONAL-ORDER}}}
{\Large {\bf{REGULATED SYSTEMS IN THE STATE SPACE}}} \\
         \vspace{10mm}
         \normalsize
{\softL{}ubom\'{\i}r DOR\v{C}\'AK, Ivo PETR\'A\v{S}, Imrich KO\v{S}TIAL} \\
               \vspace{1mm}
               {    Department of Informatics and Process Control  \\
                    BERG Faculty, Technical University of Ko\v{s}ice   \\
                    B. N\v{e}mcovej 3, 042 00 Ko\v{s}ice, Slovak Republic  \\
                    phone:      (+42195) 6025172                      \\
                    e-mail: \,{\it \{dorcak, petras, kostial \}@tuke.sk}} \\
\end{center}
\vskip 5mm
\begin{abstract}
\hspace{-6.5mm}
        In this paper we present the mathematical description and
analysis of a fractional-order regulated system in
the state space. A little historical background of our results in
the analysis and synthesis of the fractional-order dynamical
regulated systems is given.
The methods and results of simulations of the fractional-order
system described
by a state space equation equivalent to three-member
fractional-order differential equation with a fractional-order
$PD^{\delta}$
regulator are then presented.
The possibility of investigating the stability of such
systems is also considered.
\end{abstract}

\hspace{2.5mm}
{\bf Keywords:}\,
fractional calculus, fractional-order regulated system,
model, state space.

\section*{1. Introduction}

Since the first references on fractional-order derivatives in
the 17-th century, the theory of fractional-order derivatives
and integrals was  highly developed by many mathematicians.
In the last five decades, many authors made a great effort
to apply this knowledge in practice. But only in the last
decade can we find some significant works concerned with the
description, analysis, and synthesis of fractional-order
regulated systems.

      Real objects are generally fractional-order, however, for
many of them the fractionality is very low. But the control
systems used so far were all considered as integer-order systems,
regardless of the reality. Because of the higher complexity and
the absence of adequate mathematical tools, fractional-order
dynamical systems were only treated marginally in the theory and
practice of control systems, e.g. \cite{L1, Axtell, Kalojanov}.
In these works the first generalizations of analysis methods
for fractional-order
control systems were made ($s$-plane, frequency
response, etc.).

      At that time our research in this field began.
On the basis of fundamental works \cite{L1, Oldham} we made
at first
an analysis of the fractional dimension of the controlled
systems and properties of fractional regulators
\cite{Kmetek1, Prokop1, Prokop2, Kmetek2}. Until then we
did not make difference between the terms fractal and fractional
regulator.
Very significant theoretical works appeared in 1993 and later
\cite{Podlubny2, Podlubny1a, Podlubny1b, Podlubny3, Podlubny4,
Podlubny5}.
This and the following results in this field were summarized in book
\cite{Podlubny6}.
The fundamentals for practical utilization of numerical methods
\cite{Podlubny2, Podlubny1a, Podlubny1b},
based on the relation for approximation of the fractional derivatives
\cite{Rus},
and analytical methods
\cite{Podlubny3}
for fractional-order systems simulation were made.
But the proposed iterative numerical
methods did not converge to the correct solution.
The first explicit noniterative numerical methods for
simulation~~of~~fractional-order~~systems~~were presented in works
\cite{Dorcak1, Dorcak2, Dorcak3}. In works \cite{Dorcak1, Dorcak2}
the first comparison was made of these noniterative numerical
methods and corresponded analytical methods derived for
fractional-order
control system from analytical solutions published in
\cite{Podlubny3}.
An example of experimental design of a fractional-order
$PD^{\delta}$ controller, with comparison of dynamic properties
in fractional- and integer-order system with a fractional- and
integer-order controller, designed for an integer-order system
as the best approximation to given fractional-order system,
was presented in
works \cite{Dorcak1, Dorcak2}
too.
 It followed from these results, that an application of
integer-order regulator to fractional-order system is
inadequate and with a change of system or regulator parameters
can lead to system instability.
This example was taken to many other works
\cite{Podlubny4, Podlubny6, Podlubny7, Gombik, Petras0} etc.
        After experimental method of the fractional-order system
parameters identification
\cite{Dorcak1, Dorcak2},
two exact methods
\cite{Dorcak4, Podlubny8} were derived.
A great effort was devoted to elaboration of methods of
fractional-order controllers synthesis. From purely
experimental methods
\cite{Dorcak1, Dorcak2}
we continued with methods based on minimization of regulated
square
\cite{Gombik, Sykorova1}
etc., and methods based on stability and damping measure
\cite{Petras0, Petras1, Petras2} etc. In these works
the fundamentals of fractional-order
$PI^{\lambda}D^{\delta}$ controllers synthesis are developed
and now first works appear with application of these methods
to the control of chaotic fractional-order systems
\cite{Petras2, Petras3}.
         The development of the methods described in works
\cite{Sykorova1, Petras0, Petras1, Petras2}
required the elaboration of methods of fractional-order
system analysis in the frequency domain
\cite{Dorcak5, PetrasDorcak1} etc.
        The following Bode analysis was utilized for the
fractional-order controllers synthesis and such regulated
control systems stability analysis.

        This contribution is a continuation of the previous works
and deals with the mathematical description and
analysis of fractional-order regulated systems in
the state space. We present results of simulations of
a fractional-order system with the aim to investigate the stability
of such system.


\thispagestyle{empty}
\section*{2. Example~of~the~fractional-order~control~system~in~state~space}


        For the definition of the control system we consider
a simple unity feed-back control system illustrated in Fig.1,
where $G_{s}(s)$ denotes the transfer function of the controlled
system and $G_{r}(s)$ is the controller transfer function, both
integer- or fractional-order.

\vspace*{-1.9mm}

\setlength{\unitlength}{1mm}
\begin{picture}(110,42)\label{fig1}
   \put(53, 2){Fig.1: \it Feed - back control loop \rm}
   \put(26,27){$W(s) \ \ +$}
   \put(39,19){$ -$}
   \put(49.5,27){$E(s)$}
   \put(81.0,27){$U(s)$}
   \put(117.0,27){$Y(s)$}
   \put(26,25){\vector(1,0){15.9}}
   \put(45,10){\vector(0,1){11.7}}
   \put(45,25){\circle{6}}
   \put(42.9,27.0){\line(1,-1){4.0}}
   \put(42.9,23.0){\line(1, 1){4.0}}
   \put(48,25){\vector(1,0){12.2}}
   \put(45,10){\line(1,0){72.1}}
   \put(60,20){\framebox(20,10) [cc]{$G_r(s)$}}
   \put(80,25){\vector(1,0){10.0}}
   \put(90,20){\framebox(20,10) [cc]{$G_s(s)$}}
   \put(110,25){\vector(1,0){15.9}}
   \put(117,10){\line(0,1){15.1}}
\end{picture}
\vskip 4 mm
        The differential equation of the above closed regulation
system for the transfer function of the controlled system
$G_{s}(s)=1/(a_2s^{\alpha}+a_1s^{\beta}+a_0)$
and the controller
$G_{r}(s)=K+T_ds^{\delta}$
has the form
\begin{equation} \label{r1}
     a_{2}\, y^{(\alpha)}(t) +
     a_{1}\, y^{(\beta)}(t) +
     T_{d}\, y^{(\delta)}(t) +
     (a_{0}+K)\, y(t) =
     K\, w(t) +
     T_{d}\, w^{(\delta)}(t)
\end{equation}
where  $\alpha, \beta, \delta$, are generally  real
numbers and $a_0, a_1, a_2, K, T_d$ are  arbitrary
constants.

     In works
\cite{Matignon, Vinagre1, Petras2}
was presented a state space model
that expresses the fractional-order derivatives
\begin{eqnarray} \label{r2}
        {\bf x}^{(\alpha)}(t) =
           {\bf A}~{\bf x}(t) + {\bf B}~u(t), \nonumber \\
         y(t) = {\bf C}~{\bf x}(t),
   \hspace{1.5em} t \geq 0. \hspace{0.2em}
\end{eqnarray}
This description is convenient for simple models of systems
\cite{Matignon}
with only one fractional-order derivation.
%
%

 In this contribution we propose state space model of the linear
time invariant one dimensional system which expresses the first
derivatives in the state space equations (\ref{r3}) and which has
the classical state space interpretation for the fractional-order
system too. On the right side of these equations we can then
transfer more than one fractional-order derivatives
of the state space variables (${\bf x}^{(fr)}(t)$).
A disadvantage of this expression is that we cannot describe
the state space equations in vector and matrix relations as in
previous description
(\ref{r2}).
\begin{eqnarray} \label{r3}
   {\bf x}'(t) = {\bf f}({\bf x}^{(fr)}(t), u(t)),
     \hspace{1.9em}
     \nonumber \\
        y  (t) = {\bf g}({\bf x}^{(fr)}(t), u(t)),
     \hspace{1.5em} t \geq 0. \hspace{-2.1em}
\end{eqnarray}


    We verified the above methods on an example from
\cite{Dorcak1, Dorcak2}. Assume the system described by
differential equation
(\ref{r1})
with coefficients $a_2 = 0.8, a_1 = 0.5,
a_0 = 1, \alpha = 2.2, \beta = 0.9, K = 20.5, T_d = 3.7343$
and $\delta = 1.15$.
After its modification and with state space variables
$x(t) = x_1(t), x'(t) = x_2(t), x''(t) = x_2^{'}(t)$
we can derive the following state space model
equivalent to model
(\ref{r1})
\begin{eqnarray} \label{r4}
    x_1'(t) = x_2(t)~,
     \hspace{10.1cm}
     \nonumber \\
    x_2'(t) = - \frac{a_0+K}{a_2}~x_1^{(2-\alpha)}(t)
              - \frac{T_d}  {a_2}~x_1^{(1+\delta-\alpha)}(t)
              - \frac{a_1}  {a_2}~x_1^{(1+\beta-\alpha)}(t)
              - \frac{1}    {a_2}~w_1^{(2-\alpha)}(t)~,
     \hspace{0em} \\
        y  (t) =   K~x_1(t)
                 + T_d~x_2^{(\delta-1)}(t)~,
     \hspace{1.5em} t \geq 0. \hspace{17.1em}
     \nonumber
\end{eqnarray}
 We can made other alternative state space models for the same
system. The order of the integer-order system~in~the~state~space
model~was equivalent to the number of the state variables.~In
fractional-order~sys-
tems it does not hold. We propose to consider
the order of the fractional-order systems according to the order
of the highest derivative (integer- or fractional-order) in the
resulting differential equation, e.g. (\ref{r1}).

 Under condition ${\bf x}'(t)={\bf 0}$ we can obtain from
state space model
(\ref{r4})
the two equations of statics~of~this system, which are not
algebraic equations as in integer-order linear systems,but
differential equations,from which we can compute the coordinates
of the equilibrium point, to which the state trajectories
(\ref{r4})
tend.

  After discretisation of the first derivative
by first differentiation in
(\ref{r4})
we obtained the simple Euler methods for solving the state
space model. For the fractional-order derivatives in
(\ref{r4})
we can take the relation from e.g.
\cite{Podlubny2, Podlubny1a, Podlubny6, Dorcak1, Dorcak2}.
In Fig.2 and Fig.3 is a comparison of the unit - step response
of the classical numerical solution
\cite{Dorcak1, Dorcak2}
and the numerical solution of the state space model.
The obtained state trajectories represent stable focal point
for the above mentioned coefficients of the system
(with stability measure $S_t=-1.5$ and damping measure
$T_l=0.37$, e.g. \cite{Dorcak5, PetrasDorcak1})
and unstable focal point for only one changed
coefficient $T_d = 0.7343$ of the same system.

\begin{figure}\label{frequency_response}
    \vskip -5 mm
     \centerline{\psfig{file=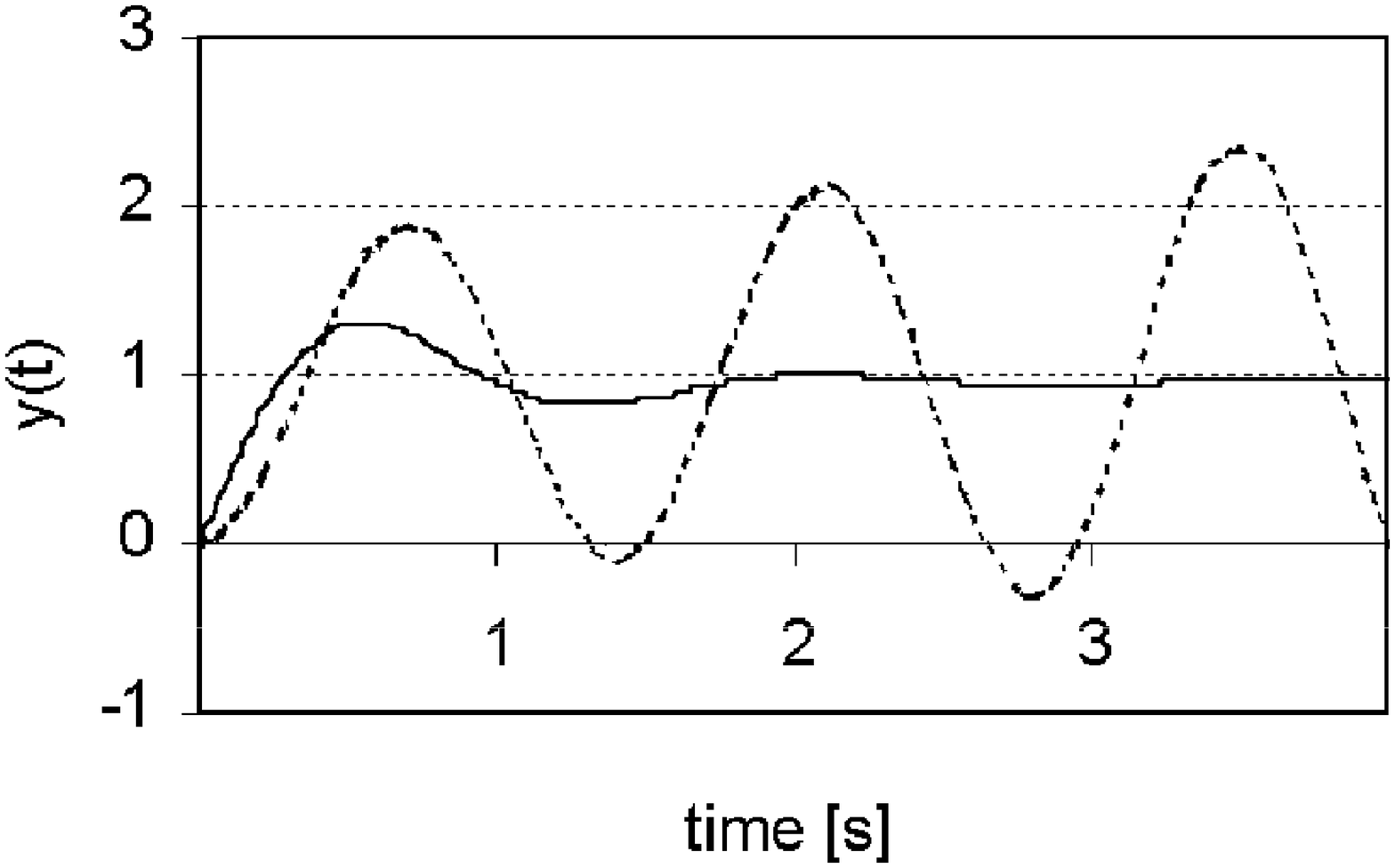,width=7.5cm}
		 \psfig{file=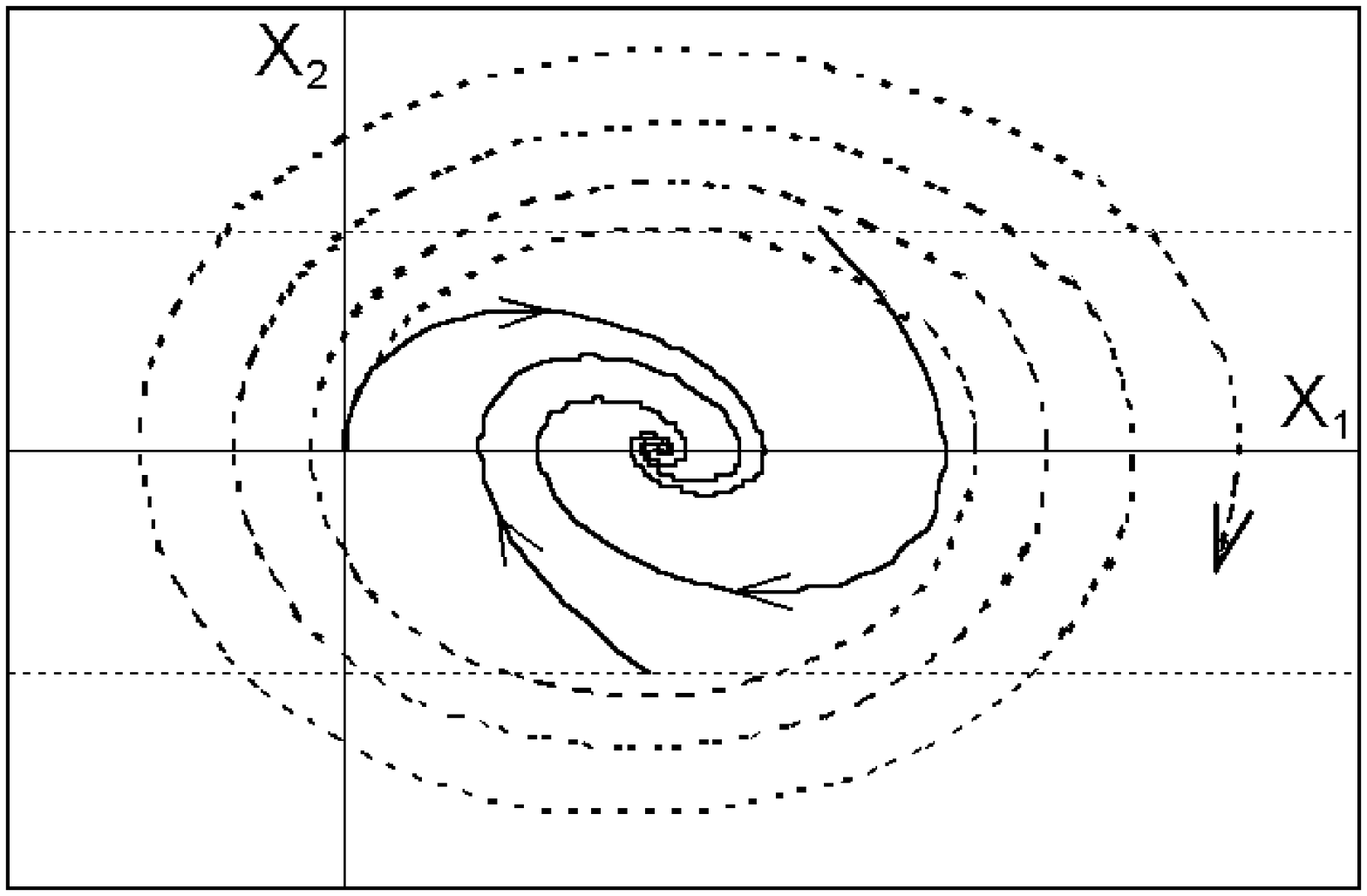,width=7.5cm}} 
    \hskip 1.5cm  \rm Fig.2 : \it Unit step responses
    \hskip 3.3cm  \rm Fig.3 : \it State trajectories \rm
    \vskip 0 mm
\end{figure}


%
%
%
%
%

This work was partially supported by grant VEGA 1/7098/20 from
the Slovak Agency for Science.

\vspace{-1.8mm}
\thispagestyle{empty}

\end{document}